\documentclass[11pt]{article}
\usepackage{mathrsfs}
\usepackage{graphicx}
\usepackage{latexsym,amsmath,amssymb,amsfonts}
\usepackage{color,colordvi,amscd,indentfirst,graphics}
\allowdisplaybreaks

\newtheorem{theorem}{Theorem}[section]

\newtheorem{example}[theorem]{Example}
\newtheorem{claim}[theorem]{Claim}

\def\qed{\hfill \rule{4pt}{7pt}}
\def\pf{\noindent {\it Proof.} }

\begin{document}

\title{Spanning trees of $K_{1,4}$-free graphs whose reducible stems have few leaves}

\author{Pham Hoang Ha\footnote{E-mail address: ha.ph@hnue.edu.vn.}\\
	Department of Mathematics\\
	Hanoi National University of Education\\
	136 XuanThuy Str., Hanoi, Vietnam\\
	\medskip\\
	Le Dinh Nam\footnote{E-mail address: nam.ledinh@hust.edu.vn}\\
	School of Applied Mathematics and Informatics\\
	Hanoi University of Science and Technology\\
	1 Dai Co Viet road, Hanoi, Vietnam\\
	\medskip\\
	Ngoc Diep Pham \footnote{E-mail address: ngocdiep23394@gmail.com (Corresponding author)}\\
	Department of Mathematics\\
	Hanoi National University of Education\\
	136 XuanThuy Str., Hanoi, Vietnam	
 }%

\date{}
\maketitle{}
\bigskip

\begin{abstract}
Let $T$ be a tree, a vertex of degree one is a \emph{leaf} of $T$
and a vertex of degree at least three is a \emph{branch vertex} of
$T$. The {\it reducible stem } of $T$ is the smallest subtree that contains all branch vertices of $T$. In this paper, we give some sharp sufficient conditions for $K_{1,4}$-free graphs to have a spanning tree whose reducible stem having few leaves.
\end{abstract}



\newpage

\section{Introduction}

In this paper, we only consider finite simple graphs. Let $G$ be a
graph with vertex set $V(G)$ and edge set $E(G)$. For any vertex
$v\in V(G)$, we use $N_G(v)$ and $\deg_G(v)$ (or $N(v)$ and $\deg (v)$ if
there is no ambiguity) to denote the set of neighbors of $v$ and the
degree of $v$ in $G$, respectively. For any $X\subseteq V(G)$, we
denote by $|X|$ the cardinality of $X$. Sometime, we denote by $|G|$ instead of $|V(G)|.$ We define
$N_G(X)=\bigcup\limits_{x\in X}N_G(x)$ and $\deg_G(X)=\sum\limits_{x\in
	X}\deg_G(x)$. We
use $G-X$ to denote the graph obtained from $G$ by deleting the
vertices in $X$ together with their incident edges. We define $G-uv$ to be the
graph obtained from $G$ by deleting the edge $uv\in E(G)$, and
$G+uv$ to be the graph obtained from $G$ by adding a new edge $uv$ joining two non-adjacent vertices $u$ and $v$ of $G$.  For two vertices $u$ and $v$ of $G$, the distance between $u$ and $v$ in $G$ is denoted 
by $d_{G}(u, v)$. We write $A:=
B$ to rename $B$ as $A$.

For an integer $m\geqslant 2,$ let $\alpha^{m}(G)$ denote the number defined by
$$\alpha^{m}(G)=\max\{ |S|:S\subseteq V(G),d_{G}(x,y)\geqslant m\,\ \text{for all distinct vertices }\,x,y\in S\}.$$
For an integer $p\geqslant 2$, we define
$$\sigma_p^{m}(G)=\min\left\lbrace \deg_G(S) : S\subseteq V(G), |S|=p, d_{G}(x,y)\geqslant m\ \text{for all distinct vertices}\;x,y\in S\right\rbrace.$$
For convenience, we define $\sigma^{m}_{p}(G)=+\infty$ if $\alpha^{m}(G)<p$. We note that,  $\alpha^{2}(G)$ is often written $\alpha(G)$, which is independence number of $G,$ and $\sigma_p^{2}(G)$ is often written $\sigma_{p}(G)$, which is minimum degree sum of $p$ independent vertices.  

Let $T$ be a tree. A vertex of degree one is a \emph{leaf} of $T$
and a vertex of degree at least three is a \emph{branch vertex} of
$T$. The set of leaves of $T$
is denoted by $L(T)$ and the set of branch vertices of $T$ is denoted by $B(T)$. The subtree $T-L(T)$ of $T$ is called the {\it stem} of $T$ and is denoted by $Stem(T)$. 

There are several sufficient conditions (such as the independence
number conditions and the degree sum conditions) for  a
graph $G$ to have  a spanning tree with a bounded number of leaves
or branch vertices. Win~\cite{Wi} obtained the following theorem, which confirms a conjecture of Las Vergnas~\cite{LV71}, and Broersma and Tunistra~\cite{BT98} gave the following sufficient condition for a graph to have a spanning tree with at most $k$ leaves. Beside that, recently, the second named author \cite{Ha1} also gave an improvement of Win by giving an independence number condition for a graph having a spanning tree which covers a certain subset of $V(G)$ and has at most $l$ leaves.

\begin{theorem}\label{t1}{\rm (Win~\cite{Wi})}
	Let $l\geqslant 1$ and $k\geqslant 2$ be integers and let $G$ be a $l$-connected graph.  If
	$\alpha(G)\leqslant k+l-1$, then $G$ has a spanning tree with at most $k$
	leaves.
\end{theorem}
\begin{theorem}\label{t2}{\rm (Broerma and Tuinstra~\cite{BT98})}
	Let $G$ be a connected graph and let $k\geqslant 2$ be an integer. If
	$\sigma_2(G)\geqslant |G|-k+1$, then $G$ has a spanning tree with at most
	$k$ leaves.
\end{theorem}

Moreover, many researchers studied spanning trees in connected graphs whose stems have a bounded number of leaves or branch vertices (see  \cite{KY}, \cite{KY15}, \cite{TZ} and \cite{Yan} for examples). 

On the other hand, for a positive integer $t \geq 3,$ a graph $G$ is said to be  $K_{1,t}-$ free graph if it contains no $K_{1,t}$ as an induced subgraph. If $t=3,$ the $K_{1,3}-$ free graph is also called the claw-free graph.   By considering the graph $G$ to be restricted in some special graph classes, many analogue researchs have been introduced (see \cite{CCHZ}, \cite{CCH14}, \cite{CHH}, \cite{HH},  \cite{KKMOSY12}, \cite{Ky09}, \cite{Ky11}  and \cite{MS84} for examples). 

For a leaf $x$ of $T$, let $y_x$ denote the  nearest branch vertex to $x$. For
every leaf $x$ of $T$, we remove the path $P_T [x, y_x)$ from $T$, where $P_T [x, y_x)$ denotes the
path connecting $x$ to $y_x$ in $T$ but not containing $y_x$. Moreover, the path $P_T [x, y_x)$ is
called the {\it leaf-branch path of $T$ incident to $x$} and denoted by $B_x$.   The resulting subtree of $T$ is called the {\it reducible stem } of $T$ and denoted by $R\_Stem(T).$ This means that the reducible stem of $T$ is the smallest subtree of $T$ which contains all branch vertices of $T$.

Recently, Ha, Hanh and Loan introduced a new concept on spanning trees. They studied the sufficient conditions for a graph to have a spanning tree whose reducible stem having a few leaves. Namely, they obtained the following theorem.
\begin{theorem}{\rm (Ha et al. \cite{HHL})}\label{thm0}
	Let $G$ be a connected graph and let $k\geqslant 2$ be an integer.  If one of the following conditions holds, then $G$ has a spanning tree with at most $k$ peripheral branch vertices. 
	\begin{enumerate}
		\item[{\rm (i)}] $\alpha(G) \leqslant 2k+2,$
		\item[{\rm (ii)}] $\sigma_{k+1}^4(G) \geqslant \left \lfloor \frac{\vert G \vert - k}{2}\right \rfloor.$
	\end{enumerate}
	Here, the notation $\lfloor r\rfloor$ stands for the biggest integer not exceed the real number $r.$
\end{theorem}
After that, Hanh also proved the following.
\begin{theorem}{\rm (Hanh \cite{Hanh})}\label{thm03}
	Let $G$ be a connected claw-free graph and let $k \geq 2$ be an integer. If one of the following conditions holds, then $G$ has a spanning tree whose reducible stem has at most $k$ leaves. 
	\begin{enumerate}
		\item[{\rm (i)}] $\alpha(G) \leq 3k+2,$
		\item[{\rm (ii)}] $\sigma_{k+1}^6(G) \geq \left \lfloor \frac{\vert G \vert -4k-2}{2}\right \rfloor.$
	\end{enumerate}
\end{theorem}
Beside that, Ha, Hanh, Loan and Pham also gave a sharp degree condition for a graph to have a spanning tree whose reducible stem having a bounded number of branch vertices.
\begin{theorem}{\rm (Ha et al. \cite{HHLP})}\label{thm0-1}
	Let $G$ be a connected graph and let $k\geqslant 2$ be an integer. If the following conditions holds, then $G$ has a spanning tree $T$ whose reducible stem has at most $k$ branch vertices. 
	\begin{equation*}
		\sigma_{k+3}^4(G) \geqslant \bigg\lfloor \frac{\vert G \vert - 2k-2}{2}\bigg\rfloor.
	\end{equation*}
\end{theorem}

In the above theorems, the authors only consider the cases of degree sum conditions with a small number of independent vertices.  A different question is whether we may give a sharp condition of $\sigma_{m}(G) (m\geq 2k+3)$ to show that $G$ has a spanning tree whose reducible stem has at most $k$ leaves. For any connected graph case, this question remains open.\\

In the case of claw-free graphs, Ha recently proved the following.
\begin{theorem}{\rm (Ha \cite{Ha2})}\label{}
	Let $G$ be a connected claw-free graph and let $k$ be an integer ($k\geq 2$).  If  $\sigma_{3k+3}(G) \geq \left\vert G \right\vert-k$, then $G$ has a spanning tree whose reducible stem has at most $k$ leaves.
\end{theorem}
In this paper, we consider the case of $K_{1,4}$-free graphs. In particular, we prove the following theorems:
\begin{theorem}\label{thm1}
	Let $G$ be a connected $K_{1,4}$- free graph.  If  $\sigma_{7}(G)) \geq \left\vert G \right\vert$, then $G$ has a spanning tree whose reducible stem is a path.
\end{theorem}
\begin{theorem}\label{thm2}
	Let $G$ be a connected $K_{1,4}$- free graph and let $k\geq 3$ be an integer.  If  $\sigma_{2k+3}(G)) \geq \left\vert G \right\vert - k+1$, then $G$ has a spanning tree whose reducible stem has at most $k$ leaves.
\end{theorem}
	To end this section, we introduce two examples to show that our results are sharp.
\begin{example}	
Let $m\geq 1$ be an integer. Let $\{R_i\}_{i=1}^3$ and $\{H_i\}_{i=1}^3$ be $6$ disjoint copies of the complete graph $K_m$ of order $m$. Let $\{x_i\}_{i=1}^{3}$ and $w$ be $4$ disjoint vertices not contained in $\cup_{i=1}^3\bigg(V(R_i)\cup V(H_i)\bigg)$. Join $x_{i}$ to all the vertices of $V(R_i)$ and all the vertices of $V(H_i)$ for every $1 \leq i \leq 3$. Adding $3$ edges $x_{1}w, x_2w, x_3w$. Let $H$ denote the resulting graph (see Figure \ref{Pic1}). Then, $H$ is a $K_{1,4}$-free graph and we have $|H|=6m+4.$ 

\begin{figure}[h]
	\centering
	\includegraphics[width=0.7\textwidth]{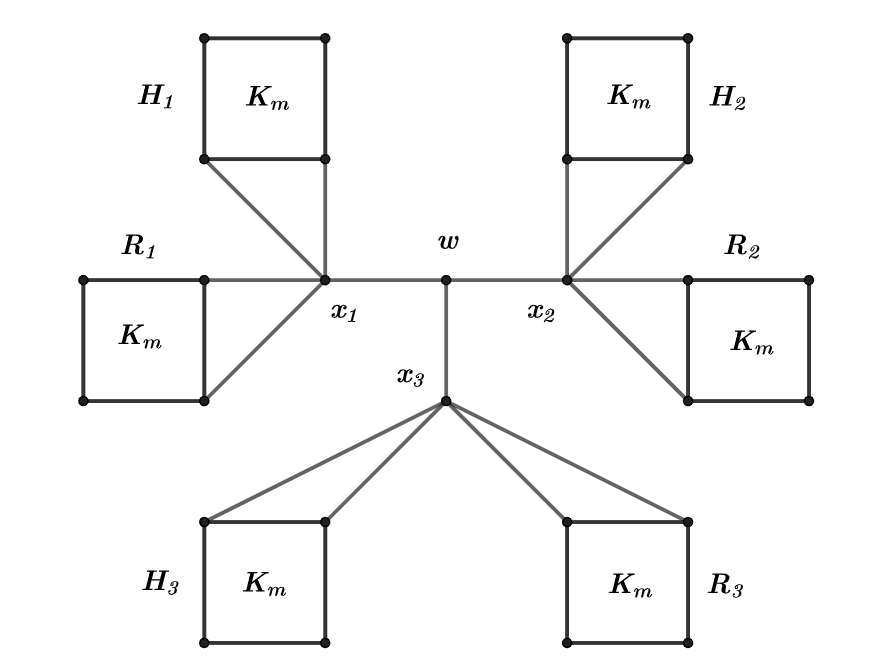}
	\caption[Graph H]{Graph $H$}
	\label{Pic1}
\end{figure}
\noindent On the other hand, take a vertex $u_i \in V(R_i)$ and a vertex $v_i\in V(H_i)$ for every $1 \leq i \leq 3.$ We obtain 
\begin{eqnarray*}
	\sigma_{7}(H)&=&\sum\limits_{i=1}^{3}\bigg(\deg_H(u_i)+\deg_H(v_i)\bigg)+\deg_H(w)\\  
	& =& 6m+3= |H|-1.
\end{eqnarray*}
But $H$ has no spanning tree whose reducible stem is a path. Hence the condition of Theorem \ref{thm1} is sharp. 
	\end{example}
 \begin{example}Let $l \geq 1$ and $m\geq 1$ be integers. Set $k=2l+1.$ Let $\{R_i\}_{i=0}^k$ and $\{H_i\}_{i=0}^{k}$ be $2k+2$ disjoint copies of the complete graph $K_m$ of order $m$. Let $D$ be a complete graph with $V(D)=\{w_i\}_{i=0}^{l}.$ Let $\{x_i\}_{i=0}^{k}$ be $k+1$ vertices not contained in $\cup_{i=0}^k\bigg(V(R_i)\cup V(H_i)\bigg)\cup \cup_{i=0}^{l}\{w_i\}$. Join $x_{i}$ to all the vertices of $V(R_i)\cup V(H_i)$ for every $0 \leq i \leq k$. Joining $w_i$ to $x_{2i}$ and $x_{2i+1}$ for all $0\leq i \leq l$. Let $G$ denote the resulting graph (see Figure 2). Then, $G$ is a $K_{1,4}$-free graph and we have $|G|=l+1+(k+1)(2m+1).$ 
\begin{figure}[h]
	\centering
	\includegraphics[width=0.7\textwidth]{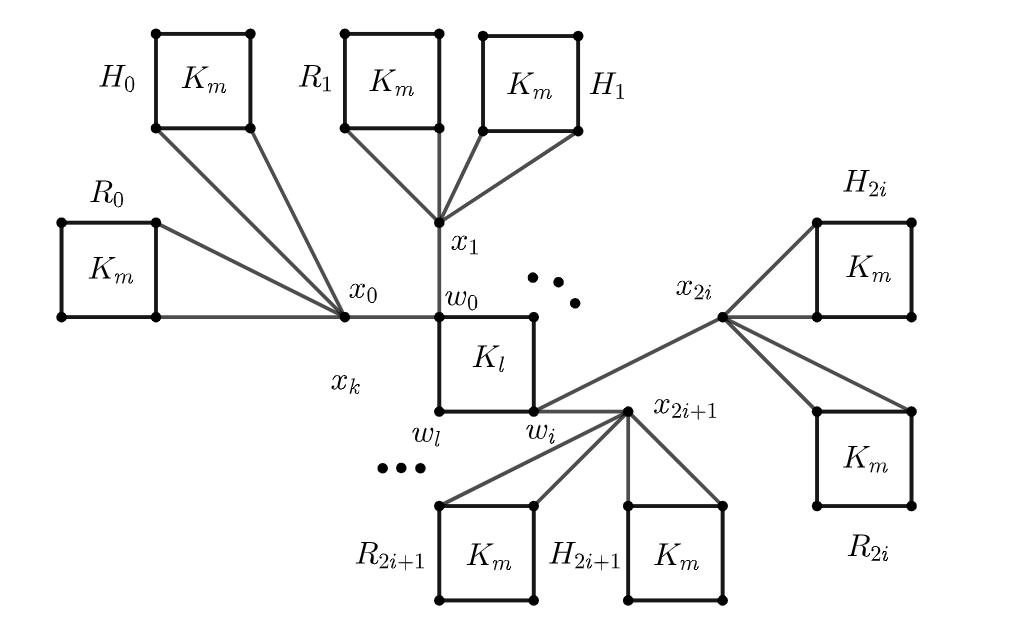}
	\caption[Graph G]{Graph $G$}
	\label{Pic2}
\end{figure}

\noindent Moreover, take a vertex $u_i \in V(R_i)$ and a vertex $v_i\in V(H_i)$ for every $0 \leq i \leq k.$ We obtain 
\begin{eqnarray*}
	\sigma_{2k+3}(G)&=&\sum\limits_{i=0}^{k}\bigg(\deg_{G}(u_i)+\deg_{G}(v_i)\bigg)+\deg_{G}(w_0)\\  
	& =&  2m(k+1)+l+2= |G|-k.
\end{eqnarray*}
But $G$ has no spanning tree whose reducible stem has at most $k$ leaves. Hence the condition of Theorem \ref{thm2} is sharp. 


\end{example}

\section{Proof of Theorems \ref{thm1} and \ref{thm2}}
 Let $T$ be a tree. For two distinct vertices $u,v$ of $T$, we denote by $P_T[u,v]$ the unique path in $T$ connecting $u$ and $v.$ We always define the
 \emph{orientation} of $P_T[u,v]$ is from $u$ to $v$. If $x \in V(P_T[u,v])$, then $x^+$ and $x^-$ denote the successor and predecessor of $x$ on $P_T[u,v]$ if they exist, respectively. For any $X\subseteq V(G),$ set $(N(X)\cap P_T[u,v])^{-} = \{x^{-} \vert x \in V(P_T[u,v])\setminus \{u\}\ \text{and}\ x \in N(X) \}$ and  $(N(X)\cap P_T[u,v])^{+} = \{x^{+} \vert x \in V(P_T[u,v])\setminus \{v\} \ \text{and}\ x \in N(X) \}.$ We refer
 to~\cite{Di05} for terminology and notation not defined here.
 
{\it Proof of Theorem \ref{thm1} and \ref{thm2}.} Suppose to the contrary that there does not exist any spanning tree $T$ of $G$ such that $\vert L(R\_Stem(T))\vert \leq k,\, (k\geq 2).$ Then every spanning tree $T$ of $G$ satisfies $\vert L(R\_Stem(T))\vert \geq k+1$.

Choose $T$ to be a spanning tree of $G$ such that
\begin{itemize}
\item [$($C0$)$]
	$\vert L(R\_Stem (T)) \vert$ is as small as possible,
\item [$($C1$)$]
$\vert R\_Stem (T) \vert$ is as small as possible subject to (C0),
\item [$($C2$)$]
  $\vert L(T)\vert $ is as small as possible subject to (C0) and (C1).
\end{itemize}

Set $L(R\_Stem(T))=\{x_1,x_2,...,x_{l}\}, l \geq k+1.$ Then, we have the following claim.
\begin{claim}\label{claim2.2}
For every $i\in \{1,2,...,l\}$, there exist at least two leaf-branch paths of $T$ which are incident to $x_i$.
\end{claim}
\begin{claim}\label{claim2.3}
	$L(R\_Stem(T))$ is an independent set in $G$.
\end{claim}
\pf
Suppose that there exist two vertices $x_i, x_j \in L(R\_Stem(T))$ such that $x_ix_j \in E(G).$ Set $H := T+x_ix_j.$ Then $H$ is subgraph of $G$ including a unique cycle $C$. Since $k\geq 2$, then $\vert L(R\_Stem (T)) \vert  \geq 3.$ Hence, we have $|B(R\_Stem(T))| \geq 1$. Hence there exists a branch vertex in the $R\_Stem(T)$ contained in $C$. Let $e$ be an edge of $C$ incident with $u$. By removing the edge $e$ we obtain a spanning tree $T'$ of $G$. Then $ \vert L(R\_Stem(T')) \vert < \vert L(R\_Stem(T)) \vert,$ the reason is that either $R\_Stem(T')$ has only one new leaf and $x_i, x_j$ are not leaves of $R\_Stem(T')$ or $x_i$ (or $x_j$) is still a leaf of $R\_Stem(T')$ but $R\_Stem(T')$ has no any new leaf and $x_j$ (or $x_i$ respectively ) is not a leaf of $R\_Stem(T')$. This contradicts the condition (C0).  The proof of Claim \ref{claim2.3} is completed.
\qed
\begin{claim}\label{claim2.4}
For each $i \in \{1,2,...,l\}$, there exist $y_i, z_i \in L(T)$ such that $B_{y_i}, B_{z_i}$ are incident to $x_i$ and $N_G(y_i) \cap (V(R\_Stem(T)) - \{x_i\}) = \emptyset$ and $N_G(z_i) \cap (V(R\_Stem(T))- \{x_i\}) = \emptyset$.
\end{claim}

\pf Let $\{a_{ij}\}_{j=1}^m$ be the  subset of $L(T)$ such that $B_{a_{ij}}$ is adjacent to $x_i$. By Claim \ref{claim2.2}, we obtain $m \geq 2.$ \\
Suppose that there are more than $m-2$ vertices $\{a_{ij}\}_{j=1}^m$ satisfying
\begin{center}
$N_G(a_{ij}) \cap \left( V(R\_Stem(T))- \{x_i\}\right) \neq \emptyset$.
\end{center}
Without loss of generality, we may assume that $N_G(a_{ij}) \cap \left( V(R\_Stem(T))- \{x_i\}\right) \neq \emptyset$ for all $j=2,...,m.$ Set $b_{ij}\in N_G(a_{ij}) \cap \left(V(R\_Stem(T))- \{x_i\}\right)$ and $v_{ij} \in N_T(x_i) \cap V(P_T[a_{ij}, x_i])$ for all $ j \in \{2,...,m\}$. Consider the spanning tree
$$T':= T + \{a_{ij}b_{ij}\}_{j=2}^m-\{x_iv_{ij}\}_{j=2}^m.$$
 Then $T'$ satisfies $|L(R\_Stem (T'))|\leq |L(R\_Stem (T))|$ and $\vert R\_Stem(T')\vert < \vert  R\_Stem(T)\vert,$ where $x_i$ is not in $V(R\_Stem(T')).$ This contradicts either the condition (C0) or (C1). Therefore, Claim \ref{claim2.4} holds.
\qed
\vskip 0.5cm
Set $U_1 =\{y_i,z_i\}_{i=1}^{l}$. 
\begin{claim}\label{claim2.5}
	$U_1$ is an independent set in $G$.
\end{claim}
\pf
Suppose that there exist two vertices $u, v \in U_1$ such that $uv \in E(G).$ Without loss of generality, we may assume that $v=y_i$ for some $i\in \{1,2,...,l\}.$ Set $v_i \in N_T(x_i)\cap V(B_{y_i}).$ Consider the spanning tree $T':= T +uy_i -v_ix_i.$ Then $|L(R\_Stem (T'))|\leq |L(R\_Stem (T))|.$ If $\deg_T(x_i)=3$ then $x_i$ is not a branch vertex of $T'.$ Hence $|R\_Stem(T')| < |R\_Stem(T)|,$ this contradicts either the condition (C0) or (C1). Otherwise, we have $|L(R\_Stem (T'))|= |L(R\_Stem (T))|$, $|R\_Stem(T')| = |R\_Stem(T)|$ and $|L(T')|< |L(T)|,$  where either $T'$ has only one new leaf $v_i$ and $y_i, u$ are not leaves of $T'$ or $y_i$ is still a leaf of $T'$ but $T'$ has no any new leaf and $u$ is not a leaf of $T'$. This contradicts the condition (C2). The proof of Claim \ref{claim2.5} is completed.
\qed
\vskip 0.5cm
Since $k\geq 2$, then $\vert L(R\_Stem (T)) \vert  \geq 3.$ Hence, we have $|B(R\_Stem(T))| \geq 1$. Let $w$ be a vertex in $B(R\_Stem(T)) .$ By Claim \ref{claim2.4} and Claim \ref{claim2.5} we conclude that $U:= U_1\cup \{w\}$ is an independent set in $G.$ This implies that $\alpha(G) \geq 2k+3.$
\vskip 0.5cm
\begin{claim}\label{claim2.6}
For every $i, j \in \{1,2,...,l\}$ where $i \neq j,$ we have $ N_G(y_i) \cap V(B_{y_j}) = \emptyset$ and $N_G(y_i) \cap V(B_{z_j}) = \emptyset$.
\end{claim}

\pf By the same role of $y_j$ and $z_j,$ we only need to prove $ N_G(y_i) \cap V(B_{y_j}) = \emptyset.$ Suppose the assertion is false. Then there exists a vertex $x \in N_G(y_i) \cap V(B_{y_j})$. Set $H := T + xy_i.$ Then $H$ is a subgraph of $G$ including a unique cycle $C,$ which contains both $x_i$ and $x_j$. 

Since $k\geq 2$, then $ |L(R\_Stem (T))| \geq 3.$ Hence, we obtain $|B(R\_Stem(T))| \geq 1$. Then there exists a branch vertex of $R\_Stem(T)$ contained in $C.$ Let $e$ be an edge incident to such a vertex in $C$ and $R\_Stem(T)$. By removing the edge $e$ from $H$ we obtain a spanning tree $T'$ of $G$ satisfying $|L(R\_Stem(T'))|<  |L(R\_Stem(T))|,$ the reason is that either $R\_Stem(T')$ has only one new leaf and $x_i, x_j$ are not leaves of $R\_Stem(T')$ or $x_i$ (or $x_j$) is still a leaf of $R\_Stem(T')$ but $R\_Stem(T')$ has no any new leaf and $x_j$ (or $x_i$ respectively ) is not a leaf of $R\_Stem(T')$. This is a contradiction with the condition (C0). So Claim \ref{claim2.6} is proved.
\qed

Now, we choose $T$ to be a spanning tree of $G$ satisfying the conditions (C0)-(C2) and the followings.
\begin{itemize}
	\item [$($C3$)$]
	$\displaystyle\sum_{i=1}^{l}\deg_T(x_i)$ is as small as possible subject to (C0)-(C2),
	\item [$($C4$)$]
	$\displaystyle\sum_{i=1}^{l}\bigg( |B_{y_i}| + |B_{z_i}| \bigg)$ is as large as possible subject to (C0)-(C3) and Claim \ref{claim2.4}, and denoted by $S_T$.	
\end{itemize}

\begin{claim}\label{claim2.7}
 For every $ p\displaystyle \in L(T)-U_1,$ we have $ \displaystyle\sum_{u\in U}|N_G(u) \cap V(B_p) | \leq \vert B_p\vert$.
\end{claim}

\pf Set $v_p \in B(T)$ to be the nearest branch vertex to $p.$ Consider the path $P_T[p,v_p],$ then we have $B_p=P_T[p,v_p^{-}].$ \\

Assume that there exists a vertex $x \in V(B_p)$ such that $xu \in E(G)$ for some $u \in U.$ Consider the spanning tree
$$T':=\left\{\begin{array}{ll}T+xu-v_p^{-}v_p,
& \;\mbox{ if } x=p,\\
T+xu -xx^+, & \;\mbox{ if } x \not= p.
\end{array}\right.$$
Then $|L(R\_Stem (T'))|= |L(R\_Stem (T))|$, $|R\_Stem(T')| = |R\_Stem(T)|,$$|L(T')|\leq |L(T)|,$ $\displaystyle\sum_{i=1}^{l}\deg_{T'}(x_i)=\displaystyle\sum_{i=1}^{l}\deg_T(x_i)$ and 	$S_{T'} > S_{T}.$ This contradicts  either the condition (C2) (if $x\in \{p, v_p^-\}$) or the condition (C4) for otherwise. This concludes that
$$
\displaystyle\sum_{u\in U}|N_G(u) \cap V(B_p) | = | N_G(w) \cap V(B_{p})|\leq |B_{p}|.
$$
Claim \ref{claim2.7} is proved.
\qed

\begin{claim}\label{claim2.8}
For every $1 \leq i \leq l$, we have $
\displaystyle\sum_{u\in U}|N_G(u) \cap V(B_{y_i}) | \leq  \vert B_{y_i} \vert- 1 $ and $\displaystyle\sum_{u\in U}|N_G(u) \cap V(B_{z_i}) | \leq  \vert B_{z_i}\vert - 1.$
\end{claim}

\pf 
As $y_i$ and $z_i$ play the same role, we only need to prove $
\displaystyle\sum_{u\in U}|N_G(u) \cap V(B_{y_i}) | \leq  \vert B_{y_i} \vert- 1 $. Set $V(B_{y_i})\cap N_T(x_i)=\{x_i^{-}\}.$ Now we consider $B_{y_i}=P_T[y_i,x_i^{-}].$ 

By Claim~\ref{claim2.6}, we obtain the following subclaim.\\
\noindent{\it Subclaim \ref{claim2.8}.1.} $ N_G(U) \cap V(B_{y_i}) =  N_G(\{y_i,z_i\}) \cap V(B_{y_i})$.\\ 

\noindent{\it Subclaim \ref{claim2.8}.2.} If $x \in N_G(y_i)\cap V(B_{y_i})$ then $x^- \notin N_G(z_i)\cap V(B_{y_i})$.

 Suppose that there exists $x \in N_G(y_i)\cap V(B_{y_i})$ such that $ x^- \in N_G(z_i)\cap V(B_{y_i})$. Consider the spanning tree $T' := T+\{xy_i,z_ix^-\}-\{xx^-,x_i^{-}x_i\}$. Then $|L(R\_Stem(T'))|\leq |L(R\_Stem(T))|.$ If $\deg_T(x_i)=3$ then $x_i$ is not a branch vertex of $T'.$ Hence $|R\_Stem(T')| < |R\_Stem(T)|,$ this contradicts either the condition (C0) or (C1). Otherwise, we have $|L(R\_Stem(T'))|= |L(R\_Stem(T))|,$ $|R\_Stem(T')| = |R\_Stem(T)|$ and $|L(T')|< |L(T)|,$  where $y_i$ and $z_i$ are not leaves of $T'.$ This is a contradiction with the condition (C2). Therefore,
Subclaim \ref{claim2.8}.2 holds.\\

\noindent{\it Subclaim \ref{claim2.8}.3.} If $x \in N_G(y_i)\cap V(B_{y_i})$ then $x^- \notin N_G(w)\cap V(B_{y_i})$.

Suppose that there exists $x \in N_G(y_i)\cap V(B_{y_i})$ such that $ x^- \in N_G(w)\cap V(B_{y_i})$. Consider the spanning tree $T' := T+\{xy_i,wx^-\}-\{xx^-,x_i^{-}x_i\}$. Then $|L(R\_Stem(T'))|\leq |L(R\_Stem(T))|.$ If $\deg_T(x_i)=3$ then $x_i$ is not a branch vertex of $T'.$ Hence $|R\_Stem(T')| < |R\_Stem(T)|,$ this contradicts either the condition (C0) or (C1). Otherwise, we have $|L(R\_Stem(T'))|= |L(R\_Stem(T))|,$ $|R\_Stem(T')| = |R\_Stem(T)|,$ $|L(T')|= |L(T)|$ and $\displaystyle\sum_{i=1}^{l}\deg_{T'}(x_i) < \displaystyle\sum_{i=1}^{l}\deg_T(x_i).$ Then the condition (C3) is false. Therefore,
Subclaim \ref{claim2.8}.3 holds.\\

\noindent{\it Subclaim \ref{claim2.8}.4.} We have $x_i^{-} \notin N_G(z_i)$ and  $x_i^{-} \notin N_G(w).$

Indeed, assume that $x_i^{-}z_i\in E(G).$ We consider the spanning tree $T':= T+x_i^{-}z_i-x_ix_i^-.$ Hence $|L(R\_Stem(T'))|\leq|L(R\_Stem(T))|,$ $ |R\_Stem(T')|\leq |R\_Stem(T)|$ and $|L(T')|<|L(T)|,$ where $z_i$ is not a leaf of $T'.$ This contradicts either the condition (C0) or (C1) or (C2).

On the other hand, if  $x_i^{-} \notin N_G(w).$ We consider the spanning tree $T':= T+x_i^{-}w-x_ix_i^-.$ Then $|L(R\_Stem(T'))|\leq |L(R\_Stem(T))|,$ $ |R\_Stem(T')|\leq |R\_Stem(T)|, |L(T')|=|L(T)|$ and $\displaystyle\sum_{i=1}^{l}\deg_{T'}(x_i) < \displaystyle\sum_{i=1}^{l}\deg_T(x_i).$ This contradicts either (C0) or (C1) or (C3). Subclaim \ref{claim2.8}.4 is proved.\\

\noindent{\it Subclaim \ref{claim2.8}.5.} If $x \in N_G(z_i)\cap V(B_{y_i})$ then $wx^+,wx^- \notin E(G)$.

Indeed, assume that $xz_i\in E(G).$ By Subclaim \ref{claim2.8}.4 and Claim \ref{claim2.5}, we obtain $x\not=x_i^{-}$ and $x\not= y_i.$ Then $x^{+}, x^{-}$ exist.\\
If $ wx^+ \in E(G)$ then we consider the spanning tree $T' := T+\{xz_i,wx^+\}-\{xx^+,x_i^{-}x_i\}$. Hence $|L(R\_Stem(T'))|\leq |L(R\_Stem(T))|.$ If $\deg_T(x_i)=3$ then $x_i$ is not a branch vertex of $T'.$ Hence $|R\_Stem(T')| < |R\_Stem(T)|,$ this contradicts either the condition (C0) or (C1). Otherwise, we have $L(R\_Stem(T'))= L(R\_Stem(T)),$ $R\_Stem(T') = R\_Stem(T),$ $|L(T')|= |L(T)|$ and $\displaystyle\sum_{i=1}^{l}\deg_{T'}(x_i) < \displaystyle\sum_{i=1}^{l}\deg_T(x_i).$ This is a contradiction with the condition (C3). Hence $wx^+ \not\in E(G).$ \\
Similarity, we also obtain $wx^- \not\in E(G).$
Therefore,
Subclaim \ref{claim2.8}.5 holds.\\

\noindent{\it Subclaim \ref{claim2.8}.6.} If $x \in N_G(z_i)\cap V(B_{y_i})$ then  $x \notin N_G(w)\cap V(B_{y_i}).$

Indeed, assume that $xw\in E(G).$ Since $xz_i \in E(G)$ and combining with Subclaim \ref{claim2.8}.4, Claim \ref{claim2.5} we obtain $x\not=x_i^{-}$ and $x\not= y_i.$ Then there are $x^{-}, x^{+}$ in $V(B_{y_i}).$ 

If $x^{-}x^{+}\in E(G),$ we consider the spanning tree $T':= T-\{x_{z_i}x_i, xx^{-}, xx^{+}\}+ \{x^{-}x^{+},xz_i,xw\}.$ Hence $|L(R\_Stem(T'))|\leq|L(R\_Stem(T))|,$ $ |R\_Stem(T')|\leq |R\_Stem(T)|,$ $|L(T')|=|L(T)|$ and $\displaystyle\sum_{i=1}^{l}\deg_{T'}(x_i) < \displaystyle\sum_{i=1}^{l}\deg_T(x_i).$ This contradicts either (C0) or (C1) or (C3). 

On the other hand, by Subclaim \ref{claim2.8}.5 we get $z_ix^{-}, z_ix^{+},wx^{-},wx^{+}\not\in E(G).$ Then $\{x^+,x^-,z_i, w\}$ is an independent set in $G$ and
$G[\{x,x^+,x^-,z_i,w\}]$ is an induced $K_{1,4}$ subgraph of $G$,
giving a contradiction. So the assertion of the claim holds. Subclaim \ref{claim2.8}.6 is proved.

By Subclaims \ref{claim2.8}.2-\ref{claim2.8}.6 we conclude that  $\{y_i\}, N_G(y_i) \cap V(B_{y_i})$ and $ \left(N_G(\{z_i,w\}) \cap V(B_{y_i})\right)^{+} $ are pairwise disjoint subsets in $B_{y_i}$. Combining with Subclaim \ref{claim2.8}.1, we have 
\begin{align*} 
\displaystyle\sum_{u\in U}|N_G(u) \cap V(B_{y_i}) |& = | N_G(y_i) \cap V(B_{y_i})|+| N_G(z_i) \cap V(B_{y_i})|+| N_G(w) \cap V(B_{y_i})|\\
&= | N_G(y_i) \cap V(B_{y_i})|+| N_G(z_i,w) \cap V(B_{y_i})|\\
&=| N_G(y_i) \cap V(B_{y_i})|+| (N_G(z_i, w) \cap V(B_{y_i}))^{+}| \leq |B_{y_i}|-1.
\end{align*}
This completes the proof of Claim~\ref{claim2.8}.
\qed
\vskip 0.5cm
{\it Proof of Theorem \ref{thm1}.} Since $G$ is $K_{1,4}-$free graph and the set $L(R\_Stem(T))$ is independent by Claim \ref{claim2.3} then $|N_G(w)\cap L(R\_Stem(T))|\leq 3.$

Combining with Claim~\ref{claim2.4} and Claims~\ref{claim2.7}-\ref{claim2.8} we obtain that 
\begin{align*}
\deg_G(U) &= \displaystyle \sum_{i=1}^{l} \left(\deg_G(y_i)+\deg_G(z_i)\right)+\deg_G(w)\\
&\leq \displaystyle \sum_{i=1}^{l}\left(\vert B_{y_i}\vert -1 \right) + \displaystyle \sum_{i=1}^{l} \left(\vert B_{z_i}\vert -1 \right)+ 2l+\displaystyle \sum_{p\in L(T)-U_1} \vert B_p \vert +\\
&+ \vert R\_Stem(T)\vert - |L(R\_Stem(T))|+3-1\\
&= \vert G \vert - \vert L(R\_Stem(T)) \vert+2\\
&\leq \vert G \vert - (k+1) +2 = \vert G \vert -k +1.
\end{align*}
When $k=2,$ we have $ \sigma_7G \leq \deg_G(U) \leq \vert G \vert - 1.$ This contradicts the assumption of Theorem \ref{thm1}. Therefore, the proof of Theorem \ref{thm1} is completed.
\qed

{\it Proof of Theorem \ref{thm2}.} 
\begin{claim}\label{claim2.9}
Let $u$ be a branch vertex of $R\_Stem(T)$  and $x_i \in L(R\_Stem(T)),i \in \{1,2,...,l\}.$ Then $x_ix \not\in E(G)$ for every $x\in N_{R\_Stem(T)}(u)\setminus V(P_{R\_Stem(T)}[x_i,u]).$
\end{claim}
\pf  If $x_ix \in E(G)$, consider the spanning tree $T' := T + x_ix-ux.$ Then $ |L(R\_Stem (T'))| < |L(R\_Stem (T))|,$ the reason is that $x_i$ is not a leaf of $R\_Stem(T').$ This contradicts the condition (C0). So Claim \ref{claim2.9} holds.
\qed

Now, we choose $T$ to be a spanning tree of $G$ satisfying the conditions (C0)-(C4) and the following
\begin{itemize}
	\item [$($C5$)$]
	$\max_{x\in B(R\_Stem(T))}\bigg\{\deg_{R\_Stem (T)}(x)\bigg\}$ is as large as possible subject to (C0)-(C4).
\end{itemize}

{\it Case 1}. There exists a vertex $w\in B(R\_Stem(T))$ such that $|N_G(w)\cap L(R\_Stem(T))|\leq 2$ then using Claim~\ref{claim2.4}, Claims~\ref{claim2.7}-~\ref{claim2.8} we obtain that 
\begin{align*}
\deg_G(U) &= \displaystyle \sum_{i=1}^{l} \left(\deg_G(y_i)+\deg_G(z_i)\right)+\deg_G(w)\\
&\leq \displaystyle \sum_{i=1}^{l}\left(\vert B_{y_i}\vert -1 \right) + \displaystyle \sum_{i=1}^{l} \left(\vert B_{z_i}\vert -1 \right)+ 2l+\displaystyle \sum_{p\in L(T)-U_1} \vert B_p \vert +\\
&+ \vert R\_Stem(T)\vert - |L(R\_Stem(T))|+2-1\\
&= \vert G \vert - \vert L(R\_Stem(T)) \vert+1\\
&\leq \vert G \vert - (k+1) +1 = \vert G \vert -k.
\end{align*} 
This contradicts the assumption of Theorem \ref{thm2}

{\it Case 2}. $|N_G(w)\cap L(R\_Stem(T))|\geq 3$ for all $w\in B(R\_Stem(T)).$ 

Since $G$ is $K_{1,4}$-free and Claim \ref{claim2.3} holds, we obtain $|N_G(w)\cap L(R\_Stem(T))|= 3$ for all $w\in B(R\_Stem(T)).$ \\
\begin{claim}\label{claim2.10}
	Let $w$ be a branch vertex of $R\_Stem(T)$ and $\{a,b,c\}=N_{R\_Stem (T)}(w)\cap L(R\_Stem(T)).$ Then there is no vertex $x\in N_{R\_Stem(T)}(w)$ such that $w\in V( P_{R\_Stem(T)}[x,a])\cap V(P_{R\_Stem(T)}[x,b])\cap V(P_{R\_Stem(T)}[x,c]).$  
\end{claim}
\pf  Assume that there exists a vertex $x\in N_{R\_Stem(T)}(w)$ such that $w\in V(P_{R\_Stem(T)}[x,a])\cap V(P_{R\_Stem(T)}[x,b])\cap V(P_{R\_Stem(T)}[x,c]).$ By Claim \ref{claim2.3} and Claim \ref{claim2.9} we obtain $\{a,b,c,u\}$ is an independent set in $G$. Hence, $G[\{w,a,b,c,x\}]$ is an induced $K_{1,4}$ subgraph of $G.$ This is a contradiction.  So Claim \ref{claim2.10} is proved.
\qed

Choose $w\in B(R\_Stem(T))$ such that $\deg_{R\_Stem (T)}(w)=\max_{x\in B(R\_Stem(T))}\bigg\{\deg_{R\_Stem (T)}(x)\bigg\}.$ 

If $\deg_{R\_Stem (T)}(w)\geq 4.$ Take $a,b,c \in L(R\_Stem(T)) \cap N_{G}(w).$ Hence there exists a vertex $x\in N_{R\_Stem(T)}(w)$ such that $w\in V(P_{R\_Stem(T)}[x,a])\cap V(P_{R\_Stem(T)}[x,b])\cap V(P_{R\_Stem(T)}[x,c]).$ This contradicts Claim \ref{claim2.10}. 

If $\deg_{R\_Stem (T)}(w)=3$ then $\deg_{R\_Stem(T)}(u) =3$ for all $u \in B(R\_Stem(T)).$ On the other hand, since $k\geq 3$, then $\vert L(R\_Stem (T)) \vert  \geq 4.$ Hence, we have $|B(R\_Stem(T))| \geq 2.$ 

\begin{claim}\label{claim2.11}
	Let $u$ be a vertex in $B(R\_Stem(T))\setminus \{w\}$ and let $a\in N_{R\_Stem (T)}(u)\setminus P_{R\_Stem(T)}[u,w].$ We have $wa\not\in E(G).$  
\end{claim}
\pf  Assume that there exists a vertex $a\in N_{R\_Stem (T)}(u)\setminus P_{R\_Stem(T)}[u,w]$ such that $wa\in E(G).$  Consider the tree $T':= T -ua+wa.$ Then $T$ satisfies the conditions (C0)-(C4) but $\deg_{R\_Stem (T)}(w)<\deg_{R\_Stem (T')}(w),$ a contradiction with the condition (C5).  So Claim \ref{claim2.11} holds.
\qed

Now, using Claim~\ref{claim2.4}, Claims~\ref{claim2.7}-\ref{claim2.8} and Claim~\ref{claim2.11} we obtain 
\begin{align*}
\deg_G(U) &= \displaystyle \sum_{i=1}^{l} \left(\deg_G(y_i)+\deg_G(z_i)\right)+\deg_G(w)\\
&\leq \displaystyle \sum_{i=1}^{l}\left(\vert B_{y_i}\vert -1 \right) + \displaystyle \sum_{i=1}^{l} \left(\vert B_{z_i}\vert -1 \right)+ 2l+\displaystyle \sum_{p\in L(T)-U_1} \vert B_p \vert +\\
&+ \vert R\_Stem(T)\vert - |L(R\_Stem(T))|+3-2\\
&= \vert G \vert - \vert L(R\_Stem(T)) \vert+1\\
&\leq \vert G \vert - (k+1) +1 = \vert G \vert -k.
\end{align*} 
This contradicts the assumption of Theorem \ref{thm2}. This completes the proof of Theorem \ref{thm2}.
\qed
\bigskip



\begin{thebibliography}{99}
	\addtolength{\baselineskip}{-1ex}
	
	\bibitem{BT98}
	Broersma, H., Tuinstra, H.: Independence trees and Hamilton cycles,
	{\it J. Graph Theory} {\bf 29} (1998),  227--237 
	
	\bibitem{CCHZ} Chen, G., Chen, Y., Hu, Z., Zhang, S.: Spanning trees with at most $k $ leaves in 2-connected $K_{1,r}$-free graphs, {\it Appl. Math. Comput.} {\bf 445} (2023), Paper No. 127842, 11 pp.
	
	
	\bibitem{CCH14}
	Chen, Y.,  Chen, G.,  Hu, Z.: Spanning $3$-ended trees in
	$k$-connected $K_{1,4}$-free graphs, {\it Sci. China Math.} {\bf 57} (2014)
	, 1579--1586. 
	
	\bibitem{CHH}
	Chen, Y., Ha, P. H., Hanh, D. D.: Spanning trees with at most 4 leaves in $K_{1,5}$-free graphs, {\it Discrete Math.} {\bf 342} (2019), 2342-2349.
	
	
	\bibitem{Di05}
	Diestel, R.: Graph Theory, 3rd Edition, Springer, Berlin, 2005.
	
	
	\bibitem{Ha1}
Ha, P. H., A note on the independence number, connectivity and $k$-ended tree, {\it Discrete Appl. Math.} \textbf{305} (2021), 142-144.

		\bibitem{Ha2}
	Ha. P. H., Spanning trees of a claw-free graph whose reducible stems have few leaves, to appear in {\it Studia Sci. Math. Hungar.} {\bf 60} (2023), no. 1, 2-15.
	
	\bibitem{HH}
	Ha, P. H., Hanh, D. D.:
	Spanning trees of connected $K_{1,t}$-free graphs whose stems have a few leaves, {\it Bull. Malays. Math. Sci. Soc.} {\bf 43}, 2373-2383 (2020)
	
		\bibitem{HHL}
	Ha, P. H., Hanh, D. D., Loan, N. T.,:
	Spanning trees with few peripheral branch vertices, {\it Taiwanese J. Math.}, Vol. {\bf 25}, No. 3, pp. 435- 447, (2021) 
	
		\bibitem{HHLP}
Ha, P. H., Hanh, D. D., Loan, N. T., Pham, N. D.:
	Spanning trees whose reducible stems have a few branch vertices, {\it Czech Math J.} {\bf 71} (146), 697-708 (2021)  
	
			\bibitem{Hanh}
	Hanh, D. D.,:
	Spanning trees with few peripheral branch vertices in a connected claw-free graph, {\it Acta Math. Hungar.} {\bf 169} (2023), no. 1, 1-14. 
		
	\bibitem{KKMOSY12}
	Kano M., Kyaw, A.,  Matsuda, H.,  Ozeki, K.,  Saito, A.,  Yamashita, T.:
	Spanning trees with a bounded number of leaves in a claw-free graph,
	{\it Ars Combin.} {\bf 103} , 137--154 (2012)
	
	\bibitem{KY} 
	Kano, M., Yan, Z.: Spanning trees whose stems have at most $k$ leaves, {\it Ars Combin.} \textbf{CXIVII}, 417-424 (2014)
	
	\bibitem{KY15} 
	Kano, M., Yan, Z.: 
	Spanning trees whose stems are spiders, {\it Graphs Combin.} \textbf{31}, 1883-1887 (2015)
	
	\bibitem{Ky09}
	Kyaw, A.:  Spanning trees with at most $3$ leaves in $K_{1,4}$-free
	graphs, {\it Discrete Math.} {\bf 309}, 6146--6148 (2009)
	
	\bibitem{Ky11}
	 Kyaw, A.: Spanning trees with at most $k$ leaves in $K_{1,4}$-free
	graphs, {\it Discrete Math.} {\bf 311}, 2135--2142  (2011)
	
	\bibitem{LV71}
	 Las Vergnas, M.: Sur une propri\'{e}t\'{e} des arbres maximaux dans un
	graphe, {\it C. R. Acad. Sci. Paris Ser. A} {\bf 272},
	1297--1300  (1971)
	
	
	\bibitem{MS84}	 Matthews, M. M.,  Sumner, D. P.: Hamiltonian results in
	$K_{1,3}$-free graphs, {\it J. Graph Theory} {\bf 8}, 
	139--146 (1984)
	
	
	
	\bibitem{TZ}
Tsugaki, M., Zhang, Y: Spanning trees whose stems have a few leaves, {\it Ars Combin.} \textbf{CXIV}, 245-256 (2014)
		

\bibitem{Yan}
Yan, Z.: Spanning trees whose stems have a bounded number of branch vertices. {\it Discuss. Math. Graph Theory} \textbf{36}, 773-778 (2016)
	
		\bibitem{Wi}
	Win, S.: On a conjecture of Las Vergnas concerning certain spanning
	trees in graphs. {\it Resultate Math.} {\bf 2},  215--224 (1979)
	
		
	
	
	
		
	
	

	

	
	
	
\end{thebibliography}
\end{document}